\documentclass[12pt, letterpaper, notitlepage]{article}
\usepackage[utf8]{inputenc}
\usepackage{caption}
\usepackage{parskip}
\usepackage{amsmath}
\usepackage{amsthm}
\usepackage{amssymb,amscd}
\usepackage[mathcal]{eucal}
\usepackage{subcaption}
\usepackage{sidecap}
\usepackage{bigints}
\usepackage{setspace}
\usepackage[all]{xy}
\usepackage{algorithm2e}
\usepackage{amsfonts}
\usepackage[square,sort,comma,numbers]{natbib}
\usepackage{hyperref}

\newtheorem{theorem}{Theorem}

\newcommand\T{\rule{0pt}{3.3ex}}       
\newcommand\B{\rule[-2.5ex]{0pt}{0pt}} 

\title{A Brief Account of\\Klein's Icosahedral Extensions}
\author{Leonardo Solanilla, Erick S. Barreto and Viviana Morales \thanks{Departamento de Matem\'aticas y Estad\'{\i}stica, Universidad del Tolima, Barrio Santa Elena, Ibagu\'e, Tolima, Colombia; leonsolc@ut.edu.co, esbarretos@ut.edu.co, vmoralesb@ut.edu.co}}
\date{February 2021}

\begin{document}

\begin{titlepage}
\maketitle

\begin{abstract}
We present an alternative relatively easy way to understand and determine the zeros of a quintic whose Galois group is isomorphic to the group of rotational symmetries of a regular icosahedron. The extensive algebraic procedures of Klein in his famous \textit{Vorlesungen über das Ikosaeder und die Auflösung der Gleichungen vom fünften Grade} are here shortened via Heymann's theory of transformations. Also, we give a complete explanation of the so-called icosahedral equation and its solution in terms of Gaussian hypergeometric functions. As an innovative element, we construct this solution by using algebraic transformations of hypergeometric series. Within this framework, we develop a practical algorithm to compute the zeros of the quintic.
\end{abstract}

\textbf{Keywords (MSC 2020):} Computational methods for problems pertaining to field theory (12-08); 	Polynomials in real and complex fields: location of zeros (algebraic theorems) (12D10); Real polynomials: location of zeros (26C10); Zeros of polynomials, rational functions, and other analytic functions of one complex variable (30C15); Numerical computation of solutions to single equations (65H05).

\end{titlepage}

\clearpage
\section{Introduction}
Let $f\in K[X]$ be an irreducible quintic with coefficients in a field $K\le\mathbb{C}$. Let $L$ be the splitting field of $f$, $K\le L\le\mathbb{C}$, and let $\textrm{Gal}(K\le L)$ denote the Galois group of $L$ over $K$. It is a well-known fact that $f$ is solvable by radicals if and only if $\textrm{Gal}(K\le L)$ is a subgroup of the Frobenius group $F_{20}$, i.e. if and only if $\textrm{Gal}(K\le L)$
is isomorphic (up to conjugacy) to $F_{20}$ of order 20, to the dihedral group $D_{10}$ of order 10, or to the cyclic group $\mathbb{Z}/5\mathbb{Z}$. Since the Galois group of certain quintics is $A_5$ (isomorphic to the group of symmetries of a regular icosahedron) and this group is not solvable, the solutions to this type of quintic cannot be expressed by field operations and radicals from the polynomial coefficients. We will say that $L$ is an icosahedral extension of $K$ if $\textrm{Gal}(K\le L)=A_5$. We assume the reader already knows the basics of Galois Theory and quintics, as they are given e.g. in Cox \cite{Cox} and the earliest works of Abel \cite{Abel} and Ruffini \cite{Ruffini}. 

In this paper we give a method to solve any quintic whose splitting field is icosahedral. Furthermore, we also provide a numerical implementation of the method. Our approach is based on Klein's \cite{Klein} foundational paper, as well as on the more recent renditions of Heymann \cite{Heymann}, Shurman \cite{Shurman} and Nash \cite{Nash}. For the numerical part, we have been inspired by Cox, Little \& O’Shea \cite{CoxLO} and Trott \cite{Trott}. Some incomplete sketchy elementary ideas on the subject can be found in \cite{AlRod} and \cite{BarMor}.
 
Tschirnhaus  transformations \cite{{Tschirnhaus}} are key ingredients of Klein's method. As usual, given a general quintic $x^{5}+c_{4}x ^{4}+c_{3}x^{3}+c_{2}x^{2}+c_{1}x+c_{0}$, the first degree transformation $z= x+ (c_{4}/5)$ eliminates the term of 4th degree and reduces the problem to a quintic $z^{5}+pz^{3}+qz^{2}+rz+s$. After a further quadratic transformation $y=z^2 - az - b$, we obtain the so-called principal quintic

\begin{equation}\label{pri_qui}
y^{5}+5\alpha y^{2}+5\beta y+\gamma=0.
\end{equation}

The quadratic transformation extends the coefficient field $K$ to $K(\delta)$, where {\small $\delta^2=\displaystyle\frac{9q^{2}}{4p^{2}}+\displaystyle\frac{3p}{5}\displaystyle \displaystyle-\frac{2r}{p}$}. For simplicity's sake, we will denote this new field by the same letter $K$. 

From a purely field-theoretical approach, our main task consists of proving the following theorem, that we will call Klein's theorem. \textit{C. f.} Slodowy \cite{Slodowy}.

\begin{theorem}
Let $K$ be a subfield of the complex numbers containing $\delta$ and the fifth roots of unity. We suppose that $L$ is a Galois extension of $K$ with $\textup{Gal}(K\le L)=A_5$, the alternating group of five letters. Then, there exists a $J\in K$ such that the solution $Y$ of the icosahedral equation $q(Y)=J$ (Section 2 below) yields $L=K(Y)$, that is, $L$ is obtained from $K$ by adjoining $Y$.
\end{theorem}

Notably, along the proof we will develop a practical numerical method for finding the zeros of our principal quintic. 

Since Klein's theorem essentially asserts that the solution of the quintic reduces to the solution of the icosahedral equation, in Section 2 we introduce this equation together with a method of solution. Gaussian hypergeometric functions provide the means to determine such a solution. Our approach is comprehensive and we provide full details. Although this solution has been described barely in Nash \cite{Nash}, we use a simpler procedure based on the algebraic transformations of the hypergeometric series. We devote Section 3 to the algebraic treatment of the principal quintic. We have closely followed the theory of Heymann's \cite{Heymann} resolvents. They allow us to show the way the solution of the quintic effectively amounts to the solution of the icosahedral equation. With this, Klein's theorem is proved. In Section 4 we rewrite the previous sections as a useful algorithm that produces as outputs the five solutions of an icosahedral quintic. In the end, we draw some conclusions regarding the main features of our method. 

\section{Solution of the icosahedral equation}
The solution of a Klein's quintic relies heavily on the geometry of the icosahedron.

\subsection*{Icosahedral symmetry}\label{icosol}
The icosahedron is the surface of the regular polyhedron or Platonic solid having 12 vertices, 30 edges and 20 triangular faces (the ancient Greek number $\epsilon\iota\kappa o\sigma\iota$, \textit{eíkosi}, means twenty). Once it is inscribed in a sphere, we project radially from the center to obtain the image of the icosahedron on the sphere. Then, we project the circumsphere  onto the complex plane via the stereographic projection to get a plane image of the icosahedral surface. As it is usual, the extended complex plane or Riemann sphere $\mathbb{S}$ is conveniently identified with the complex projective space $\mathbb{P}^1$.

The icosahedral rotation group $\mathcal{I}$ has order 60 and is isomorphic to $A_5$, the group of the even permutations of five letters. It comprises all spherical transformations or rotations under which the icosahedron remains invariant. By virtue of the identification $\mathbb{S}\simeq\mathbb{P}^1$, the elements of $\mathcal{I}$ can be understood as Möbius transformations or, alternatively, as homogeneous forms. As well, $A_5$ is isomorphic to the projective special linear group $\textrm{PSL}(2,5)$. The action of the icosahedral group on the sphere yields three abnormal orbits corresponding to the vertices (12 points), edge midpoints (30 points) and face centers (20 points). The rest of the orbits are normal in the sense that the group acts transitively on them and so, they all consist of 60 points. The icosahedral equation arises from the problem of inverting the branched covering $\mathbb{P}^1\rightarrow \mathbb{P}^1/\mathcal{I}$. 

\subsection*{Invariant forms of 60th degree}
Klein (1884) describes an inverse of this covering in relation to polynomials or homogeneous forms. For instance, the vertices of the icosahedron are precisely the zeros of the invariant polynomial
\begin{equation}\label{ico_fff}
f=f(z,w)= zw(z^{10}+11z^{5}w^{5}-w^{10}).
\end{equation}
Similarly, the Hessian determinant of $f$ produces the invariant form
\begin{equation}\label{ico_hhh}
H= H(z,w)=-( z^{20}+w^{20})+228(z^{15}w^{5}-z^{5}w^{15})-494z^{10}w^{10},
\end{equation} 
which vanishes at the centers of the icosahedron faces. Also, the Jacobian determinant of $f,H$ provides the 30th degree invariant
\begin{equation}\label{ico_ttt}
T=T(z,w)= (z^{30}+w^{30})+ 522(z^{25}w^{5}-z^{5}w^{25})-10005(z^{20}w^{10}+z^{10}w^{20}),
\end{equation}
vanishing at the midpoints of the edges. We should standardize these invariant forms in the 60th degree in order to deal with the spherical orbits under the action of $\mathcal{I}$. Once we do so, the set $\mathbb{O}$ of the invariant forms of 60th degree is easily furnished with a structure of vector space of dimension two over $\mathbb{C}$. Among other choices, the forms $H^3$ and $f^5$ form a basis for $\mathbb{O}$. Remarkably enough, $T^2=1728f^5-H^3$ in this basis. 

It is not hard to establish an identification between the homogeneous polynomials in $\mathbb{O}$ with the orbits in $\mathbb{P}^1/\mathcal{I}$, where they vanish.

\subsection*{The isomorphism $\mathbb{P}^1\simeq \mathbb{P}^1/\mathcal{I}$}
The canonical projection $\mathbb{P}^1\rightarrow \mathbb{P}^1/\mathcal{I}$ takes a point $(z,w)$ to its orbit $[z,w]$. It corresponds to a map $\mathbb{P}^1\rightarrow \mathbb{O}$, taking $(z,w)$ to the form $P(z,w)$, which equals zero at $(z,w)$. An arbitrary $P\in\mathbb{O}$ is written $P=\mu H^3-\lambda f^5$, in the basis $H^3, f^5$. Then, $P=0$ implies
$$\frac{H^3}{f^5}=\frac{\lambda}{\mu}.$$ 
So, the covering can be written $(z,w)\mapsto (\lambda,\mu)$. In non-homogeneous coordinates we have {\small $\left(Y=\frac{z}{w},1\right)\mapsto \left(\frac{\lambda}{\mu}=J,1\right)$}, or simply $Y\mapsto J$. 

The inversion of the covering reveals interesting facts. First, we must restrict ourselves to one of the possible branches or restricted domains. Then, we construct the inverse map $[z,w]\mapsto (z,w)$ by means of $\left(H^3(z,w),f^5(z,w)\right)\mapsto (z,w)$ and thus $J\mapsto Y$. In the subset of the normal orbits $[z,w]\in \mathbb{P}^1/\mathcal{I}$, it is possible to choose a unique pair of holomorphic functions $(z,w)$. However, the  abnormal orbits introduce discontinuities corresponding to the vertices, edge midpoints and face midpoints. The nature of this inverse map is next elucidated with the aid of Complex Analysis.

\subsection*{Differential resolvents}
Even more, we can construct an explicit inverse map $J\mapsto Y$ of the covering. For sure, $H^3$ and $f^5$ are connected through variable $J$ and the constant rank theorem reduces the problem to the equivalent form 
\begin{eqnarray*}
f(z,w) &=& k,\\
H(z,w) &=& u,
\end{eqnarray*}
where $k$ is constant and $u$ is a complex variable. We recall that the jacobian determinant of the map $\mathbb{P}^{1}\rightarrow \mathbb{P}^{1}$, $(z,w)\mapsto(k,u)$, equals $-20T$, where $T$ is the icosahedral form (\ref{ico_ttt}). Also, the derivatives of $z,w$ with respect to $u$ can be computed implicitly through
$$
\left( \begin{array}{l}
\frac{\displaystyle dz}{\displaystyle du} \\ 
\frac{\displaystyle dw}{\displaystyle du}  
\end{array} \right)
=-\frac{1}{20T}
\left( \begin{array}{rr}
 \frac{\displaystyle \partial H}{\displaystyle \partial w} &-\frac{\displaystyle \partial f}{\displaystyle \partial w} \\ 
-\frac{\displaystyle \partial H}{\displaystyle \partial z} & \frac{\displaystyle \partial f}{\displaystyle \partial z} 
\end{array} \right)
\left( \begin{array}{l}
\displaystyle 0 \\
\displaystyle 1
\end{array} \right).
$$
Since
$$\frac{d}{du}\left(\frac{\partial f}{\partial z}\right)=\frac{\partial^{2} f}{\partial z^{2}}\cdot \frac{dz}{du}+\frac{\partial^{2} f}{\partial z\partial w}\cdot \frac{dw}{du}\quad\textrm{and}\quad\frac{d}{du}\left(\frac{\partial f}{\partial w}\right)=\frac{\partial^{2} f}{\partial w^{2}}\cdot \frac{dw}{du}+\frac{\partial^{2} f}{\partial z\partial w} \cdot \frac{dz}{du},$$ this yields the linear ordinary homogeneous equations
\begin{equation}\label{dif_equ1}
T^{2}\frac{d^{2}z}{du^{2}}+T\frac{dT}{du}\cdot \frac{dz}{du}+\frac{11Hz}{400}=T^{2} \frac{d^{2}w}{du^{2}}+T\frac{dT}{du}\cdot \frac{dw}{du}+\frac{11Hw}{400}=0.
\end{equation}
By using $T^{2}=12^{3}f^{5}-H^{3}$ and setting $u=H$, we obtain 
\begin{equation}\label{dif_equ2}
(12^{3}f^{5}-H^{3}) \frac{d^{2}z}{dH^{2}}-\frac{3}{2}H^{2}\frac{dz}{dH}+\frac{11Hz}{400}=0
\end{equation}
and a similar equation for $w$. Then, we suitably recover $J=H^3/12^3f^5$ and perform the change of variables $H\mapsto J$ to get
\begin{equation*}
\frac{dz}{dH}=\frac{3J}{H}\cdot\frac{dz}{dJ}\quad \textrm{and} \quad
\frac{d^{2}z}{dH^{2}}=\frac{9J^{2}}{H^{2}}\cdot\frac{d^{2}z}{dJ^{2}}+\frac{2}{H}\cdot \frac{dz}{dH}.
\end{equation*}
The substitution of these derivatives in (\ref{dif_equ2}) produces the well-known hypergeometric equation
{\small
\begin{equation}\label{dif_equ3}
J(1-J)\frac{d^{2}z}{dJ^{2}}+[c-(a+b+1)J]\frac{dz}{dJ}-abz = 0,\quad\textrm{with}\quad a=\frac{11}{60},\ b=-\frac{1}{60},\ c=\frac{2}{3},
\end{equation}
}
and a corresponding equation for $w$. This equation has three regular singularities located respectively at $J=0, 1, \infty$. They correspond to the branch points of the covering: $J=0$ means $H=0$, $J=1$ means $T=0$ and $J=\infty$ means $f=0$. A general exposition of the nature of the solutions to general hypergeometric equations similar to (\ref{dif_equ3}) near ordinary points and regular singularities can be found in Copson \citep[Chapter X]{Copson} and Hille \citep[Section 6.1]{Hille}. 

The solutions to (\ref{dif_equ3}) can be expressed in terms of the hypergeometric series
\begin{equation*}\label{ghser}
F(a,b,c;J)=1+\sum_{n=1}^{\infty}\frac{(a)_n(b)_n}{(1)_n(c)_n}J^n,
\end{equation*}
where $(p)_n=p(p+1)\cdots(p+n-1)$. For the values of $a,b,c$ in (\ref{dif_equ3}), there are fundamental systems $z,w$ of solutions at the singular points $J=0,1,\infty$. In particular, at $J=\infty$, since $a-b$ is not an integer, 
$$
z(J) = \frac{1}{J^{\frac{11}{60}}} F\left(\frac{11}{60},\frac{31}{60},\frac{6}{5},\frac{1}{J}\right),\quad
w(J) = J^{\frac{1}{60}} F\left(-\frac{1}{60},\frac{19}{60},\frac{4}{5},\frac{1}{J}\right)
$$
are linearly independent for $|J|>1$. It is a known fact that the functions $z,w$ have an analytical continuation in the whole plane, save for the singularities.

\subsection*{Solving for $Y$}
Inasmuch as our main interest lies in the ratio of functions $z,w$ and not in the actual solutions to the differential equation, we consider the function
\begin{equation}\label{s}
s = s(J) = \frac{z(J)}{w(J)}=\frac{F\left(\frac{11}{60},\frac{31}{60},\frac{6}{5},\frac{1}{J}\right)}{J^{\frac{1}{5}}~F\left(-\frac{1}{60},\frac{19}{60},\frac{4}{5},\frac{1}{J}\right)}.
\end{equation}
With it, the solution to the icosahedral equation is accomplished by composing $s$ with a Möbius transformation
$$\mathfrak{M}(s)=\mathfrak{M}(s(J))=Y(J)=Y.$$ 
The spherical symmetry $\mathfrak{M}$ rotates $s$ to the correct value of $Y$, which satisfies the correct ``initial or boundary conditions'' at three convenient points. We can determine $\mathfrak{M}$ assisted by Table \ref{threeps}.

\begin{table}[h!]
\centering
 \begin{tabular}{p{5em}|c|c|c|} 
 $J$  & $0$ & $1$  & $\infty$\T\B \\  
 \hline
 $s=s(J)$ & $\infty$ & $-(-1)^{\frac{1}{10}}\sqrt[5]{1728}$ & $0$  \T\B \\
 \hline
 $Y=\mathfrak{M}(s)$ & $\infty$ & $-(-1)^{\frac{1}{10}}$ & $0$   \T\B\\ 
\hline
\end{tabular}
\caption{Determining $Y$.}
\label{threeps}
\end{table}

The first row contains three different values of $J$. The second row is calculated from the first row by using the formula defining $s(J)$ and some hypergeometric identities that we will state shortly. The values of $Y$ in the third row are obtained from $J=H^3/1728f^5$. For instance, $J=\infty$ implies $f=0$ and this is achieved through $f(Y_{\infty}=0,1)=0$.  For the case $J=0$ we realize that, at the infinity, $$\frac{1}{J}=\frac{H^3}{1728f^5}(1,Y_0=\infty)=\infty.$$
Finally, an easy calculation proves that 
$$J=\frac{H^3}{1728f^5}\left({\displaystyle -(-1)^{\frac{1}{10}}},1 \right)=1.$$
The values of $s(\infty), s(0)$ are easily determined. $s(1)$ is found with the aid of the following algebraic transformations of the Gaussian hypergeometric function, \textit{cf.} Vidunas \citep[Section 6.3]{Raimundas}:  
\begin{align}\label{iden7}
F\left(\frac{11}{60},\frac{31}{60},\frac{6}{5},\varphi_{1}(x)\right)&=\frac{(1-228x+494x^{2}+228x^{3}+x^{4})^{\frac{11}{20}}}{1+11x-x^{2}},
\end{align}
\begin{align}\label{iden8}
F\left(-\frac{1}{60},\frac{19}{60},\frac{4}{5},\varphi_{1}(x)\right)&=(1-228x+494x^{2}+228x^{3}+x^{4})^{-\frac{1}{20}},
\end{align}
where 
\begin{equation}\label{fi}
\varphi_{1}(x)=\frac{1728x(x^{2}-11x-1)^{5}}{(1-228x+494x^{2}+228x^{3}+x^{4})^{3}}.
\end{equation}
The quotient of (\ref{iden7}) into (\ref{iden8}) provides a feasible expression for (\ref{s}):
\begin{equation}\label{iden10}
s(J)=\frac{F\left(\frac{11}{60},\frac{31}{60},\frac{6}{5},\varphi_{1}(x)\right)}{J^{\frac{1}{5}}~F\left(-\frac{1}{60},\frac{19}{60},\frac{4}{5},\varphi_{1}(x)\right)}= \frac{(1-228x+494x^{2}+228x^{3}+x^{4})^{\frac{3}{5}}}{J^{\frac{1}{5}}(1+11x-x^{2})}.
\end{equation}
Thus, we may take at convenience $x = i$ to obtain $\varphi_{1}(i)=1=1/J$, that is $J=1$. Therefore,
\begin{eqnarray*}
s(1) =\frac{F\left(\frac{11}{60},\frac{31}{60},\frac{6}{5},\varphi_{1}(i)\right)}{(1)^{\frac{1}{5}}~F\left(-\frac{1}{60},\frac{19}{60},\frac{4}{5},\varphi_{1}(i)\right)} &=& \frac{(1-228(i)+494(i)^{2}+228(i)^{3}+(i)^{4})^{\frac{3}{5}}}{(1)^{\frac{1}{5}}(1+11(i)-(i)^{2})}\\
&=&-(-1)^{\frac{1}{10}}\sqrt[5]{1728}.
\end{eqnarray*}
Finally, we write 
$$Y = \mathfrak{M}(s)=\frac{\mathfrak{a}s+\mathfrak{b}}{\mathfrak{c}s+\mathfrak{d}},$$
for some $\mathfrak{a},\mathfrak{b},\mathfrak{c},\mathfrak{d}\in\mathbb{C}$. Since $s(\infty)=0$ yields $Y=0$, we must have $\mathfrak{b}=0$ and $\mathfrak{d}\neq 0$. Without loss of generality, we set $\mathfrak{a}=1$. Then, $s(0)=\infty$ and $Y=\infty$ imply $1/\mathfrak{c}=\infty$ and so, $\mathfrak{c}=0$. This leaves us with $s(1)/\mathfrak{d}=-(-1)^{\frac{1}{10}}$ and an elementary calculation gives $\mathfrak{d}=\sqrt[5]{1728}$. Hence, $Y$ is given by
\begin{equation}
Y = \mathfrak{M}(s) =\frac{1}{\sqrt[5]{1728}}\frac{z(J)}{w(J)}=\frac{F\left(\frac{11}{60},\frac{31}{60},\frac{6}{5},\frac{1}{J}\right)}{\sqrt[5]{1728\cdot J} \times F\left(-\frac{1}{60},\frac{19}{60},\frac{4}{5},\frac{1}{J}\right)}.
\end{equation}
This non-homogeneous coordinate solves the problem. 

\section{Proof of Klein's Theorem}

\subsection*{Heymann's resolvents}
The coefficients $\alpha, \beta, \gamma$ of the principal quintic $y^{5}+5\alpha y^{2}+5\beta y+\gamma=0$ belong to a quadratic Tschirnhaus extension $K$ of the original coefficient field. We remind that 
$$\Delta = 3125\cdot(108\alpha^{5}\gamma - 135\alpha^{4}\beta^{2} + 90\alpha^{2}\beta\gamma^{2} - 320\alpha\beta^{3}\gamma + 256\beta^5 + \gamma^{4})$$
is the so-called discriminant of the quintic. By hypothesis, $K$ also contains the fifth roots of unity and $\sqrt{\Delta}$, because $\textrm{Gal}(K\le L) = A_{5}$.

First and foremost, the immemorial symmetric relation $yz=y+z$ transforms this quintic into
{\small
$$
(1+5\alpha +5\beta +\gamma)z^{5}-5(3\alpha +4\beta +\gamma)z^{4} +5(3\alpha +6\beta +2\gamma)z^{3}-5(\alpha +4\beta +2\gamma)z^{2}+5(\beta +\gamma)z-\gamma.
$$
}
The fourth and third degree terms vanish whenever $\alpha=\gamma/3$ and $\beta=-\gamma/2$ (the actual solution does not have to fulfill these strident conditions). By renaming $y=\eta_1$, $z=\eta_2$, $h_1=-6\gamma^{-1}$ and $h_2=1-h_1$, we obtain Heymann's $\eta$-resolvents
\begin{eqnarray*}\label{hey_eta}
h_{1}\eta_{1}^{5}-10\eta_{1}^{2}+15\eta_{1}-6 = 0, \\
h_{2}\eta_{2}^{5}-10\eta_{2}^{2}+15\eta_{2}-6 = 0.
\end{eqnarray*}
We observe that $h_1+h_2=1$ and $\eta_1\eta_2=\eta_1+\eta_2$. Many higher degree polynomials in $\eta_1,\eta_2$ are easily reduced to quadratic or linear polynomials by these identities. For example, $h_{1}\eta_{1}^{4}\eta_{2} = 10\eta_{1}+\eta_{2}-6$. Such polynomials are called Heymann's simultaneous resolvents. 

\subsection*{The first form of the solution}
We are interested in solutions of the form
\begin{equation}\label{fir_sol}
y = p\eta_1 + q\eta_2, \quad p,q\in\mathbb{C}.
\end{equation}
The substitution of this expression in (\ref{pri_qui}) yields the simultaneous resolvent
$$M\eta_{1}^{2}+N\eta_{2}^{2}+P\eta_{1}+Q\eta_{2}+R,$$ 
where $M, N, P, Q$ and $R$ are given by
\begin{eqnarray*}
M &=& 5p^{2}(2h_{2}p^{3}+2h_{1}q^{3}+\alpha h_{1}h_{2}), \\
N &=&5q^{2}(2h_{1}q^{3}+2h_{2}p^{3}+\alpha h_{1}h_{2}),\\
P &=& 5p(-3h_{2}p^{4}+10h_{2}p^{3}q+6h_{1}pq^{3}+h_{1}q^{4}+2\alpha h_{1}h_{2}q+\beta h_{1}h_{2}), \\
Q &=& 5q(-3h_{1}q^{4}+10h_{1}q^{3}p+6h_{2}qp^{3}+h_{2}p^{4}+2\alpha h_{1}h_{2}p+\beta h_{1}h_{2}), \\
R &=& 6(h_{2}p^{5}-5h_{2}p^{4}q+10h_{2}q^{2}p^{3}+10h_{1}p^{2}q^{3}-5h_{1}pq^{4}+ h_{1}q^{5})+\gamma h_{1}h_{2}.
\end{eqnarray*}
In order to find values of the coefficients $\alpha, \beta, \gamma$, we require $M=N=R=0$. With this,
\begin{eqnarray*}\label{abg_pqh}
\alpha &=& -2(p^{3}h_{1}^{-1}+q^{3}h_{2}^{-1}), \notag \\
\beta &=& 3(p^{3}h_{1}^{-1}(p-2q)-q^{3}h_{2}^{-1}(2p-q)),\\
\gamma  &=& -6(p^{3}h_{1}^{-1}(p^{2}-5pq+10q^{2})+q^{3}h_{2}^{-1}(10p^{2}-5pq+q^{2})).\notag
\end{eqnarray*}
These formulas provide a map $(p,q,h_1)\mapsto (\alpha,\beta,\gamma)$.  We also need certain map $(\alpha,\beta,\gamma)\mapsto (p,q,h_1)$. To do so, the previous equations for $\alpha$ and $\beta$ produce
\begin{equation*}\label{hhf_pqa}
h_1 = -\frac{18p^{3}(p-q)}{3\alpha(2p-q)-2\beta}, \quad
h_2 = -\frac{18q^{3}(p-q)}{3\alpha(p-2q)+2\beta}.
\end{equation*}
By replacing these expressions in the equation for $\gamma$,
\begin{eqnarray}\label{rsh_abg}
12\alpha r+6\beta s-\gamma=0,
\end{eqnarray}
where we have introduced the new variables $r=(p-q)^2$, $s=p+q$ (this variable $s$ is different from variable $s(J)$ in Section 2). Clearly, $r$ and $s$ depend on each other. Now we replace $h_1, h_2$ in the fundamental relation $h_1+h_2=1$ and use (\ref{rsh_abg}) to obtain the central quadratic relation
\begin{equation}\label{qua_equ}
(\alpha^{4}+\alpha \beta \gamma-\beta ^{3})(12r)^{2}-(2\alpha ^{3}\gamma +11\alpha^{2}\beta^{2}+\beta \gamma^{2})(12r)+(\alpha \gamma -8\beta ^{2})^{2}=0.
\end{equation}
We notice that the quadratic discriminant needed to solve for $r$ is 
$$108\alpha^{5}\gamma-135\alpha^{4}\beta^{2}+90\alpha^{2}\beta\gamma^{2}-320\alpha\beta^{3}\gamma+256\beta^{5}+\gamma^{4}=\frac{\Delta}{3125}.$$
Therefore, the solution of (\ref{qua_equ}) does not extend our field $K$. Each value of $r$ provides a value of $s$ and so, values of $p$ and $q$. Heymann's parameters $h_1, h_2$ are the solutions of the product-sum or trinomial-factoring relations $h_1+h_2=1$ and
\begin{equation}\label{pro_hhh}
h_{1}h_{2}=\frac{432\beta p^{3}q^{3}}{12(\alpha \gamma -\beta^{2})r-\gamma^{2}}.
\end{equation}

\subsection*{Eulerian resolvents and icosahedral invariants forms}
Heymann's ideas recall some results of Euler (1764) on solvable quintics of the form (\ref{pri_qui}). His solutions have the form
$$y=\epsilon z-\epsilon^2 w,$$
where $\epsilon$ denotes a fifth root of unity. In particular, $\epsilon=1$ yields $y= z- w$. Also, $z,w$  are related to the polynomial coefficients through 
\begin{equation*}\label{eul_rel}
\alpha = -zw^{2},\quad \beta = z^{3}w,\quad \gamma=-(z^{5}-w^{5}).
\end{equation*}
Euler's solution is not our solution. However, we can use the map $(\alpha,\beta,\gamma)\mapsto (p,q,h_1)$ and give it a shot. First, we notice that the second degree coefficient  of equation (\ref{qua_equ}) vanishes with these expressions of $\alpha, \beta, \gamma$. The remaining linear equation gives 
\begin{equation*}
r = \frac{w^4(7z^{5}+w^{5})^{2}}{12f},
\end{equation*}
where $f$ is the well-known icosahedral invariant form (\ref{ico_fff}). From (\ref{rsh_abg}) it is easy to find the corresponding value
$$s=\frac{z^{3}(-z^{10}+39z^{5}w^{5}+26w^{10})}{6f}.$$
With $r$ and $s$ we determine $p$ and $q$. Then, from (\ref{pro_hhh}) we get
\begin{equation}\label{ico_par}
4 h_{1}h_{2}=\frac{H^{3}}{12^{3}f^{5}},
\end{equation}
where $H$ is the expression in (\ref{ico_hhh}). In other words, $J=4h_1h_2$ is the icosahedral parameter. 

Now, from $h_1+h_2=1$ and (\ref{ico_par}) we obtain
\begin{equation}\label{dif_hhh}
h_{1}-h_{2} =\frac{\sqrt{12^{3}f^{5}-H^{3}}}{24f^{2}\sqrt{3f}}=\frac{T}{24f^{2}\sqrt{3f}},
\end{equation}
where $T$ is just but the invariant (\ref{ico_ttt}).

The smooth way to compute the parameters $h_1, h_2$ is through $r$ and $s$:
\begin{align*}\label{R}
h_{1}=-\frac{9(s+\sqrt{r})^{3}\sqrt{r}}{6\alpha(s+3\sqrt{r})-8\beta},~~~~~~~~~~~~~
h_{2}= \frac{9(s-\sqrt{r})^{3}\sqrt{r}}{6\alpha(s-3\sqrt{r})-8\beta}.
\end{align*}
Henceforth, by virtue of (\ref{dif_hhh}),
\begin{equation}\label{dif_hhi}
h_{1}-h_{2}=\kappa(r)\sqrt{r}=\frac{T}{24f^{2}\sqrt{3f}},
\end{equation}
where $\kappa=\kappa(r)$ is a rational function of $r$ and so,
$$\sqrt{3rf}=\frac{T}{24 \kappa f}.$$ 
Since $h_1-h_2\in K$, neither $\sqrt{r}$ nor $\sqrt{3f}$ introduces extra elements into field $K$. 

\subsection*{The second form of the solution}
The main idea to find the solutions of the quintic runs as follows. From coefficients $\alpha,\beta,\gamma$ of (\ref{pri_qui}), we find $r, s$. With them, we determine $p,q$ together with $h_1,h_2$. Then, we pose the problem by forcing the solution to have also Euler's form but not its original meaning. That is,
\begin{equation}\label{sys_hey}
y = p \eta_{1}+q \eta_{2}=z-w,\quad \eta_{1}\eta_{2}= \eta_{1}+\eta_{2}.
\end{equation}
The order pair $(z,w)$ represents homogeneous coordinates of $\mathbb{P}^{1}$ and should be found. The resolution for Heymann's resolvents gives
\begin{equation*}\label{eta_1y2}
\eta_{1}=-\frac{2\sqrt{3f}}{t-\sqrt{3f}} \quad \textrm{and} \quad
\eta_{2}=\frac{2\sqrt{3f}}{t+\sqrt{3f}},
\end{equation*} 
where $t$ is the well-konwn octahedral form 
$$t(z,w)=z^{6}+2z^{5}w-5z^{4}w^{2}-5z^{2}w^{4}-2zw^{5}+w^{6}.$$
In this way, the quadratic term $\eta_1\eta_2$ in (\ref{sys_hey}) does not introduce new elements to field $K$. Further, the  solution $y$ to the reduced quintic is obtained by substituting $\eta_1, \eta_2$ in system (\ref{sys_hey}), thus
\begin{equation*}\label{mag_for}
y=-\displaystyle \frac{6sf+2t\sqrt{3rf}}{t^{2}-3f}.
\end{equation*}
We observe that $r, s$ are fixed by the coefficients of the quintic and $f$ is invariant under the action of the icosahedral group $\mathcal{I}$. The only quantity that changes is $t$. Certainly, the action of $\mathcal{I}$ yields the five values
\begin{equation}\label{mag_for1}
t_{\nu} =\epsilon^{3\nu}z^6+2\epsilon^{2\nu}z^5w-5\epsilon^{\nu}z^4w^2 -5\epsilon^{4\nu}z^2w^4-2\epsilon^{3\nu}zw^5+\epsilon^{2\nu}w^6,
\end{equation}
in which $\epsilon=e^{2\pi i/5}\in K$ and $\nu=0,1,2,3,4$. In brief, the solution set is
\begin{equation}\label{mag_for2}
\left\{ y_{\nu}=-\displaystyle \frac{6sf+2t_{\nu}\sqrt{3rf}}{t_{\nu}^{2}-3f} : \nu=0,1,2,3,4 \right\}.
\end{equation}
Now, from $\alpha,\beta,\gamma$ and the values $p,q$, we determine the icosahedral parameter
\begin{equation}\label{ico_parj}
J=4 h_{1}h_{2}=\frac{H^{3}(z,w)}{12^{3}f^{5}(z,w)}.
\end{equation}
With this, the problem is reduced to determine a pair $(z,w)$ from an orbit $[z,w]\in \mathbb{P}^{1}/\mathcal{I}$, i.e. a solution of the icosahedral equation. The corresponding $Y=Y(J)$ is determined as in Section \ref{icosol}. Then,  we replace $f(Y,1), H(Y,1)$ and $t_{\nu}(Y,1)$ in (\ref{mag_for2}) to find the zeros of the principal quintic. 

In such a way, Klein's theorem has been proved. 

\section{Algorithm}
The foregoing procedure forms a finite sequence of computer-implementable instructions.

\begin{algorithm}
\SetAlgoLined
\SetKwInOut{Input}{Read}
\Input{coefficients $\alpha,\beta,\gamma$ of a reduced quintic with Galois group $A_5$;}
\BlankLine
\SetKwFunction{finhom}{f}
\SetKwProg{Fn}{Function}{:}{}
\Fn{\finhom{$z$}}{
     \KwData{$z$, a complex number}
     \KwResult{inhomogeneous icosahedral form $f(z)=f(z,1)$}
     $f(z)\leftarrow z(z^{10}+11z^{5}-1)$\;
     \KwRet\;
}
\BlankLine
\SetKwFunction{tinhom}{t}
\SetKwProg{Fn}{Function}{:}{}
\Fn{\tinhom{$z$}}{
	 \textbf{input:} $z$, a complex number\\
	 \textbf{output:} inhomogeneous octahedral forms $t_{\nu}(z)=t_{\nu}(z,1)$, a vector of five components\\
     $\epsilon\leftarrow \exp\left(\frac{2\pi i}{5}\right)$\;
	 \For{$\nu\leftarrow 0$ \KwTo $4$}{
	 $t_{\nu}(z)\leftarrow \epsilon^{3\nu}z^6+2\epsilon^{2\nu}z^5-5\epsilon^{\nu}z^4 -5\epsilon^{4\nu}z^2-2\epsilon^{3\nu}z+\epsilon^{2\nu}$\;
	 }
	 \KwRet\;
}
\BlankLine
\SetKwFunction{funhyp}{F}
\SetKwProg{Fn}{Function}{:}{}
\Fn{\funhyp{$a,b,c;z$}}{
     \KwData{argument $z$; parameters $a,b,c$}
     \KwResult{Gaussian hypergeometric function $F(a,b,c;z)=_2F_1(a,b,c;z)$}
     $F(a,b,c;z)\leftarrow \displaystyle 1+\sum_{n=1}^{\infty}\frac{(a)_n(b)_n}{(1)_n(c)_n}z^n$\;
     \KwRet\;
}
\BlankLine
\tcc{initialization}
\BlankLine
$r\leftarrow $roots of $(\alpha^{4}+\alpha \beta \gamma-\beta ^{3})(12r)^{2}-(2\alpha ^{3}\gamma +11\alpha^{2}\beta^{2}+\beta \gamma^{2})(12r)+(\alpha \gamma -8\beta ^{2})^{2}=0$\;
$s\leftarrow $solutions of $12\alpha r+6\beta s-\gamma=0$\;
$p\leftarrow \frac{1}{2} \left(\sqrt{r}+s\right)$; $q\leftarrow \frac{1}{2}  \left(s-\sqrt{r}\right)$\;
$J\leftarrow \displaystyle 4 h_{1}h_{2}=4\frac{432\beta p^{3}q^{3}}{12(\alpha \gamma -\beta^{2})r-\gamma^{2}}$\;
$Y\leftarrow \displaystyle \frac{F\left(\frac{11}{60},\frac{31}{60},\frac{6}{5},\frac{1}{J}\right)}{\sqrt[5]{1728\cdot J} \times F\left(-\frac{1}{60},\frac{19}{60},\frac{4}{5},\frac{1}{J}\right)}$\;
\BlankLine
\tcc{computation of zeros}
\BlankLine
\For{$\nu\leftarrow 0$ \KwTo $4$}
{$y_{\nu}\leftarrow \displaystyle -\frac{6sf(Y)+2t_{\nu}(Y)\sqrt{3rf(Y)}}{t_{\nu}^{2}(Y)-3f(Y)}$\;}
\BlankLine
\SetKwInOut{Output}{Write}
\Output{zeros $y_{\nu}$, $\nu=0,1,2,3,4$, of the quintic.}
\BlankLine
\end{algorithm}
\clearpage

\subsection*{Example}
We consider the following quintic with coefficients in $\mathbb{Q}[i, \sqrt{\Delta}]$:
\begin{equation*}
y^5+5 i y^2-12 y+(1-i)=0.
\end{equation*}
We know its Galois group is $ A_5 $. We have $\alpha=i, \beta=-\frac{12}{5}, \gamma=1-i$. Firstly, we conveniently pick the solution $r=-0.140712-1.06363 i$ of
\begin{equation*}
(12 r)^2 \left(\alpha^4+\alpha \beta \gamma-\beta^3\right)-(12 r) \left(2 \alpha^3 \gamma+11 \alpha^2 \beta^2+\beta \gamma^2\right)+\left(\alpha \gamma-8 \beta^2\right)^2=0.
\end{equation*}
The corresponding value $s =0.816914 -0.0478157 i$ is obtained from $12\alpha r+6\beta s-\gamma=0$. Next, $p=\frac{1}{2} \left(\sqrt{r}+s\right)$ and $q=\frac{1}{2}  \left(s-\sqrt{r}\right)$ provide
$$~~~~~p=0.749812\, -0.413396 i~~~~~~\textrm{and}~~~~~~q=0.0671022\, +0.365581 i.~~~~$$
With all this, we obtain immediately 
$$J=4h_{1}h_{2}=\frac{H^3}{1728 f^5}=-0.324158-2.04659 i.$$
A solution of this icosahedral equation is
\begin{eqnarray*}
Y=Y(-0.324158 - 2.04659i) &=& \frac{_{2}F_{1}(\frac{11}{60},\frac{31}{60} ;\frac{6}{5}; J^{-1})}{\sqrt[5]{ 1728 J} ~ _{2}F_{1}(-\frac{1}{60}, \frac{19}{60};\frac{4}{5}; J^{-1})}\\&=&0.178352\, +0.0718131 i.
\end{eqnarray*}
Afterwards, we calculate the values $ f(Y) $ and $ t_{\nu}(Y) $:
\begin{align*}
f&=-0.178721-0.0713975 i,           &        t_{0}&= 0.509555\, -0.278001 i,\\ 
 t_{1}&= -0.761539+0.997924 i,  & t_{2}&=0.372515\, -1.14707 i,  \\
  t_{3}&=0.240993\, +1.27692 i,   &      t_{4}&=-0.361523-0.849771 i.
\end{align*}
Finally, we substitute these values in
$$y_{\nu}=-\frac{6s f(Y) + 2 t_{\nu}(Y) \sqrt{3r f(Y) }}{t_{\nu}^2(Y)-3 f(Y)}$$
to get the desired solutions
\begin{align*}
y_{0}&= 0.0895118\, -0.0828539 i,   &  y_{1}&= -0.0120031+2.20094 i,\\
y_{2}&=-0.0430531-1.43083 i,  &   y_{3}&=-1.90456-0.333135 i,\\
y_{4}&=1.87011\, -0.354121 i.
\end{align*}
This example has been developed with the routines of the program \textsl{Wolfram Ma\-the\-ma\-ti\-ca}, version 12.0.0.0.

\section{Concluding remarks}
We have constructed a simple algorithm for the evaluation of the zeros of any quintic polynomial whose splitting field is icosahedral. The mathematical machinery used to construct the method is not that arduous and more abridged than the usual approaches employed to compute such zeros. In particular, the theory of Heymann \cite{Heymann} permits to rapidly establish a central relation between the solutions to the equation and the icosahedral invariants. Also, the algebraic transformations of hypergeometric series furnish an alternative simplifying way to solve the icosahedral equation. 

Some other interesting questions arise from these results. They will be reported at a later date.

\section*{Acknowledgments}
This research was partially funded by the \textit{Comit\'e Central de Investigaciones, Universidad del Tolima, Ibagu\'e, Colombia}, grant number 60120. We also thank the \textit{Facultad de Ciencias, Universidad del Tolima}, for its logistic support to craft the manuscript. 

\vspace{6pt}  


\textbf{Conflicts of interest:} The authors declare no conflict of interest. The founding sponsors had no role in the design of the study, in the analyses, in the writing of the manuscript and in the decision to publish the results.


\end{document}